\documentclass[11pt]{article}
\usepackage{amsmath,amsthm,epsfig,amssymb}
\usepackage{amssymb}
\usepackage{latexsym}
\usepackage{graphicx}
\usepackage{epsf}

\usepackage[all]{xy}
\CompileMatrices

\newtheorem{theorem}{Theorem}[section]
\newtheorem{proposition}[theorem]{Proposition}

\newtheorem{definition}[theorem]{Definition}
\newtheorem{corollary}[theorem]{Corollary}
\newtheorem{remark}[theorem]{Remark}
\newtheorem{example}[theorem]{Example}

\newcommand{\dis}{\displaystyle}
\font\fat = msbm10 at 11pt

\def\bbb#1{\hbox {{\fat #1}}}

\def\R{{\bbb R}}

\def\Z{{\bbb Z}}
\def\S{{\bbb S}}

\input xy
\xyoption{all}

\date{\today}
\title{\bf Fundamental Groups of Commuting Elements in Lie Groups}
\author{\bf Enrique Torres-Giese\thanks{Supported by a Conacyt 
fellowship}\ \ \ and\ \ Denis Sjerve\thanks{Research supported by NSERC Discovery grant A7218}}

\begin{document}
\maketitle

\begin{abstract} We compute the fundamental group of the spaces of ordered commuting $n$-tuples
of elements in the Lie groups $SU(2)$, $U(2)$ and $SO(3)$. For $SO(3)$ the mod-2 cohomology of 
the components of these spaces is also obtained.
\end{abstract}

\section {\bf Introduction}\label{intro}
In this paper we calculate the fundamental groups of the connected components of the spaces  
$$M_n(G):=Hom(\Z^n,G),\ \mbox{ where $G$ is one of $SO(3),\ SU(2)$ or $U(2)$.}$$ 
The space $M_n(G)$ is just the space of ordered commuting $n$-tuples of elements from $G,$  topologized as a subset of $G^n.$

The spaces  $M_n(SU(2))$ and $M_n(U(2))$ are connected (see~\cite{AC}), but $M_n(SO(3))$ has 
many components if $n>1.$  One of the components is the one containing the element  $(id,id,\ldots,id);$ see Section~\ref{thespaces}. The other components are all homeomorphic to $V_2(\mathbb R^3)/\mathbb Z_2\oplus\mathbb Z_2,$ where $V_2(\mathbb R^3)$ is the Stiefel manifold of orthonormal $2$-frames in $\mathbb R^3$ and the action of 
$\mathbb Z_2\oplus\mathbb Z_2$ on  $V_2(\mathbb R^3)$ is given by
$$(\epsilon_1,\epsilon_2)(v_1,v_2)=(\epsilon_1 v_1,\epsilon_2 v_2),\ 
\mbox{where $\epsilon_j=\pm 1$ and $(v_1,v_2)\in V_2(\mathbb R^3). $}
$$
Let $e_1,e_2,e_3$  be the standard basis of $\mathbb R^3.$  Under the homeomorphism $SO(3)\to V_2(\mathbb R^3)$ given by $A\mapsto (Ae_1,Ae_2)$ the action of
$\mathbb Z_2\oplus\mathbb Z_2$ on $V_2(\mathbb R^3)$ corresponds to the action defined by right multiplication by the elements
of the group generated by the transformations 
$$(x_1,x_2,x_3)\mapsto(x_1,-x_2,-x_3),\ (x_1,x_2,x_3)\mapsto(-x_1,x_2,-x_3).$$ 
The orbit space of this action is homeomorphic to $\S^3/Q_8$, where $Q_8$ is 
the  quaternion group of order eight.

Then $M_n(SO(3))$ will be a disjoint union of many copies of  $\S^3/Q_8$ and the component containing $(id,\ldots,id).$ For brevity let $\vec{1}$ denote the $n$-tuple $(id,\ldots,id).$ 

\begin{definition}\label{comps}
{\em Let $M_n^+(SO(3))$ be the component of $\vec{1}$ in $M_n(SO(3))$, and let $M_n^-(SO(3))$ be the complement
$M_n(SO(3))-M_n^+(SO(3))$.}
\end{definition}
 
Our main result is the following

\begin{theorem}\label{mainth}
For all $n\ge 1$
$$\begin{array}{rcc} 
\pi_1(M_n^+(SO(3))) & = & \Z_2^n \\
\pi_1(M_n(SU(2))) & =& 0 \\
\pi_1(M_n(U(2))) & = & \Z^n \end{array}$$
\end{theorem}
The other components of $M_n(SO(3)),$ $n>1,$ all have fundamental group $Q_8.$

\begin{remark}\label{cover}{\em To prove this theorem we first prove that 
$\pi_1(M_n^+(SO(3))) =  \Z_2^n,$ and then use the following property of  spaces of homomorphisms (see~\cite{G}).  
Let $\Gamma$ be a discrete group, $p:\tilde{G}\to G$  a covering 
of Lie groups, and $C$ a component of the image of the induced map
$p_*:Hom(\Gamma,\tilde{G})\to Hom(\Gamma,G)$.  Then $p_*:p_*^{-1}(C)\to C$ is a regular covering with
covering group $Hom(\Gamma,Ker\ p)$.
Applying this to the universal coverings $SU(2)\to SO(3)$ and $SU(2)\times\R \to U(2)$ induces coverings
$$\Z_2^n\to M_n(SU(2))\to M_n^+(SO(3))$$ $$\Z^n\to M_n(SU(2))\times\R^n\to M_n(U(2))$$}
\end{remark}

\begin{remark}{\em The spaces of homomorphisms arise in different contexts (see~\cite{J}). 
In physics for instance, the orbit space $Hom(\Z^n,G)/G$, with $G$ acting by conjugation, is the moduli space of isomorphism 
classes of flat connections on principal $G$-bundles over the $n$-dimensional torus. Note that, if $G$ is connected,
the map $\pi_0(Hom(\Z^n,G))\to\pi_0(Hom(\Z^n,G)/G)$ is a bijection of sets. The study of these spaces arises 
from questions concerning the understanding of the structure of the components of this moduli space and their number. 
These questions are part of the study of the quantum field theory of gauge theories over the $n$-dimensional
torus (see~\cite{BFM},\cite{KS}).
}\end{remark}

The organization of this paper is as follows. In Section 2 we study the toplogy of $M_n(SO(3))$ and
compute its number of components. In Section 3 we prove Theorem~\ref{mainth} and apply this result to 
mapping spaces. Section 4 treats the cohomology of $M_n^+(SO(3))$. Part of the content of this paper is part 
of the Doctoral Dissertation of the first author (\cite{E}).

\section{\bf The Spaces $M_n(SO(3))$}\label{thespaces}

In this section we describe the topolgy of the spaces $M_n(SO(3)),\ n\ge 2.$ 
If $A_1,A_2$ are commuting elements from $SO(3)$ then there are $2$ possibilities:
either $A_1,A_2$ are rotations about a common axis; or
$A_1,A_2$ are involutions about axes meeting at right angles.
 The first possibility covers the case where one of $A_1,A_2$ is the identity since the identity can be considered as a rotation about any axis.  

It follows that there are 2 possibilities for an $n$-tuple  
$(A_1,\ldots,A_n)\in M_n(SO(3)):$ 
\begin{enumerate}
\item Either $A_1,\ldots,A_n$ are all rotations about a common axis $L$; or
\item There exists at least one pair $i,j$ such that $A_i,A_j$ are involutions about perpendicular axes.  If $v_i,v_j$ are unit vectors representing these axes then  all the other $A_k$ must be one of $id,A_i,A_j$ or $A_iA_j=A_jA_i$ (the involution about the cross product $v_i\times v_j).$
\end{enumerate}

It is clear that if $\omega(t)=(A_1(t),\ldots,A_n(t))$ is a path in 
$M_n(SO(3))$ then exactly one of the following 2 possibilities occurs:
either the rotations $A_1(t),\ldots,A_n(t)$ have a common axis $L(t)$ for all $t$; or there exists a pair $i,j$ such that  $A_i(t),A_j(t)$ are involutions about perpendicular axes for all $t$. In the second case the pair $i,j$ does not depend on $t.$

\begin{proposition}\label{compvec{1}}
$M_n^+(SO(3))$ is the space of $n$-tuples $(A_1,\ldots,A_n)\in SO(3)^n$ such that all the $A_j$ have a common axis of rotation.
\end{proposition}
\begin{proof}
Let $A_1,\ldots,A_n$  have a common axis of rotation $L$.  Thus 
$A_1,\ldots,A_n$ are rotations about $L$ by some angles  $\theta_1,\ldots,\theta_n.$  We can change all angles to $0$ by a path (while  maintaining the common axis). Conversely, if $\omega(t)=(A_1(t),\ldots,A_n(t))$ is a path containing $\vec{1}$ then  the $A_j(t)$ will have a common axis of rotation for all $t$ (which might change with $t$). 
\end{proof}

Thus any  component of  $M_n^-(SO(3))$ can be represented by an $n$-tuple $(A_1,\ldots,A_n)$ satisfying possibility 2 above.

\begin{corollary} The connected components of $M_2(SO(3))$ are $M_2^{\pm}(SO(3)).$ \end{corollary} 
\begin{proof} Let $(A_1,A_2)$ be a pair in $M_2^-(SO(3)).$ Then there are unit vectors $v_1,v_2$ in $\R^3$ such that $v_1,v_2$ are
perpendicular and $A_1,A_2$ are involutions about $v_1,v_2$ respectively. The pair $(v_1,v_2)$ is not unique since any one of the 
four pairs  $(\pm v_1,\pm v_2)$ will determine the same involutions. In fact there is a 1-1 correspondence between pairs $(A_1,A_2)$ 
in $M_2^-(SO(3))$ and sets $\{(\pm v_1,\pm v_2) \}.$  Thus $M_2^-(SO(3))$ is homeomorphic to the orbit space 
$V_2(\R^3)/\Z_2\oplus\Z_2.$  Since $V_2(\R^3)$ is connected so is $M_2^-(SO(3)).$ 
\end {proof}
 
Next we determine the number of components of $M_n^-(SO(3))$ for $n>2$. The following example will give an indication  of the complexity.

\begin{example}\label{example1}
Let $(A_1,A_2,A_3)$ be an element of $M_3^-(SO(3)).$ Then there exists at least one pair $A_i,A_j$ without a 
common axis of rotation. For example suppose $A_1,A_2$ don't have a common axis. Then 
$A_1,A_2$ are involutions about perpendicular axes $v_1,v_2$, and there are 
$4$ possibilities for $A_3:$  
$A_3=id,A_1,A_2\ \mbox{or}\ A_3=A_1A_2.$
 We will see that the   triples 
$$(A_1,A_2,id),(A_1,A_2,A_1),(A_1,A_2,A_2),(A_1,A_2,A_1A_2)$$ belong to different components. Similarly if $A_1,A_3$ or $A_2,A_3$ don't have a common axis of rotation.  This leads to 12  components, but some of them are the same component.
An analysis  leads to a total of 7 distinct components corresponding to the following 7 triples:
$(A_1,A_2,id)$, $(A_1,A_2,A_1)$, $(A_1,A_2,A_2)$, $(A_1,A_2,A_1A_2)$, $(A_1,id,A_3)$, $(A_1,A_1,A_3)$, $(id,A_2,A_3)$;
where $A_1,A_2$ are distinct involutions in the first 4 cases;  $A_1,A_3$ are distinct involutions in the next 2 cases; and $A_2,A_3$ are
distinct involutions in the last case. These components are all homeomorphic to $\S^3/ Q_8$. 
Thus $M_3(SO(3))$ is homeomorphic to the disjoint union  
$$M_3^+(SO(3))\sqcup \S^3/Q_8\sqcup\ldots\sqcup \S^3/Q_8,$$ where there are 7 copies of $\S^3/Q_8$.
\end{example}

The pattern of this example holds for all $n\ge 3.$ A simple analysis shows that 
 $M_n^-(SO(3))$ consists of $n$-tuples 
$\vec{A}=(A_1,\ldots,A_n)\in SO(3)^n$ satisfying the following conditions:
\begin{enumerate}
\item Each $A_i$ is either an involution about some axis $v_i,$ or the 
identity.
\item If $A_i,A_j$ are distinct involutions then their axes are at right 
angles.
\item There exists at least one pair $A_i,A_j$ of distinct involutions.
\item If $A_i,A_j$ are distinct involutions then every other $A_k$ is one of $id,A_i,A_j$ or $A_iA_j.$
\end{enumerate}
This leads to 5 possibilities for any element 
$(B_1,\ldots,B_n)\in M_n^-(SO(3)):$
$$
(B_1,B_2,*,\ldots,*),(B_1,id,*,\ldots,*),(id,B_2,*,\ldots,*),
(B_1,B_1,*,\ldots,*),(id,id,*,\ldots,*),
$$

where $B_1,B_2$ are distinct involutions about perpendicular axes and the asterisks are choices from amongst  $id,B_1,B_2,B_3=B_1B_2.$ The choices must satisfy the conditions above.  

These 5 cases account for all components of  $M_n^-(SO(3)),$ but not all choices lead to distinct components. If
$\omega(t)= (B_1(t),B_2(t),\ldots,B_n(t))$ is a path in 
$M_n^-(SO(3))$ then it is easy to verify the following statements:
\begin{enumerate}
\item If some $B_i(0)=id$ then $B_i(t)=id$ for all $t$.
\item If $B_i(0)=B_j(0)$ then $B_i(t)=B_j(t)$ for all $t$.
\item If $B_i(0),B_j(0)$ are distinct involutions then so are $B_i(t),B_j(t)$ for all $t$.
\item If $B_k(0)=B_i(0)B_j(0)$ then $B_k(t)=B_i(t)B_j(t)$ for all $t$.
\end{enumerate}
These 4 statements are used repeatedly in the proof of the next theorem.
\begin{theorem} 
The number of components of $M_n^-(SO(3))$ is
$$\left\{\begin{array}{ll}
\dis\frac{1}{6}(4^{n}-3\times 2^{n}+2) & \mbox{if $n$ is even}\\
 & \\
\dis \frac{2}{3}(4^{n-1}-1)-2^{n-1}+1   &\mbox{if $n$ is odd}
\end{array}\right.
$$
Moreover, each component is homeomorphic to  $\S^3/Q_8.$
\end{theorem}
\begin{proof}
Let $x_n$ denote the number of components.  The first 3 values of $x_n$ are  $x_1=0,$ $ x_2=1$ and $x_3=7,$ in agreement with the statement in the theorem. We consider the above 5 possibilities one by one. First assume $\vec{B}=(B_1,B_2,*,\ldots,*).$ Then different choices of the asterisks lead to different components.  Thus the contribution in this case is $4^{n-2}.$ Next assume 
$\vec{B}=(B_1,id,*,\ldots,*).$ Then all choices for the asterisks are admissible, except for those choices involving only $id$ and $B_1.$ This leads to 
$4^{n-2}-2^{n-2}$ possibilities.  However, changing every occurrence of $B_2$ to $B_3,$ and $B_3$ to $B_2,$ keeps us in the same component.  Thus the total contribution in this case is $(4^{n-2}-2^{n-2})/2.$ This is the same contribution for cases 3 and 4.
Finally, there are $x_{n-2}$ components associated to  $\vec{B}=(id,id,*,\ldots,*).$ 
This leads to the recurrence relation 
$$\dis x_n= 4^{n-2}+\frac{3}{2}(4^{n-2}-2^{n-2})+x_{n-2}$$ 
Now we solve this recurrence relation for the $x_n.$

Given any  element $(B_1,\ldots,B_n)\in M_n^-(SO(n))$ we select a pair of
involutions, say $B_i,B_j,$ with perpendicular axes $v_i,v_j.$ All the other $B_k$ are determined uniquely by $B_i,B_j.$ Thus the element  $(v_i,v_j)\in V_2(\mathbb R^3)$ determines $(B_1,\ldots,B_n).$ But all the elements  $(\pm v_i,\pm v_j)$ also determine $(B_1,\ldots,B_n).$ Thus the component to which $(B_1,\ldots,B_n)$ belongs is homeomorphic to 
$V_2(\mathbb R^3)/\mathbb Z_2\oplus\mathbb Z_2\cong \S^3/Q_{8}.$
\end{proof}

\section{\bf Fundamental Group of $M_n(G)$}  

In this section we prove Theorem \ref{mainth}, and we start by finding an appropriate description of $M_n^+(SO(3))$. 
Let $T^n=(\S^1)^n$ denote the $n$-torus. Then

\begin{theorem}\label{quotient}
$M_n^+(SO(3))$ is homeomorphic to the quotient space $\S^2\times T^n/\sim,$ where $\sim$ is the equivalence relation generated by
$$\displaystyle
(v,z_1,\dots,z_n)\sim (-v,\bar{z}_1,\ldots,\bar{z}_n)\ \mbox{ and }\ 
(v,\vec{1})\sim(v^{\prime},\vec{1})\ 
\mbox{for all $v,v^{\prime}\in \S^2,z_i\in \S^1.$}
$$
\end{theorem}
\begin{proof}
If $(A_1,\ldots,A_n)\in M_n^+(SO(3))$ then there exists $v\in \S^2$ such that $A_1,\ldots,A_n$ are rotations about $v$.  Let $z_j\in \S^1$ be the elements corresponding to these rotations.  The $(n+1)$-tuple $(v,z_1,\ldots,z_n)$ is not unique.  For example, if one of the $A_i$'s is not the identity then $(-v,\bar{z}_1,\ldots,\bar{z}_n)$ determines the same $n$-tuple of rotations.  On the other hand, if all the $A_i$'s are the identity then any element $v\in \S^2$ is an axis of rotation. 
\end{proof}

We will use the notation $[v,z_1,\ldots,z_n]$ to denote the equivalence class of $(v,z_1,\dots,z_n).$ 
Thus $x_0=[v,\vec{1}]\in \S^2\times T^n/\sim$ is a single point, which we choose to be the base point.  It corresponds to the $n$-tuple $(id,\dots, id)\in M_n^+(SO(3))$. Then
$M_n^+(SO(3))$ is locally homeomorphic to $\mathbb R^{n+2}$ everywhere except at the point $x_0$ where it is singular.

{\it Proof of Theorem 1.2}:  Notice that the result holds for $n=1$ since $Hom(\Z,G)$ is homeomorphic to $G$. 
The first step is to compute $\pi_1(M_n^+(SO(3))).$  Let  $T^n_0=T^n-\{\vec{1}\}$ and 
$M_n^+=M_n^+(SO(3)).$
Removing the singular point $x_0=[v,\vec{1}]$ from $M_n^+$  we have  
$M_n^+-\{x_0\}\cong \S^2\times T_0^n/\mathbb Z_2,$ see Theorem \ref{quotient}.
If $t$ denotes the generator of $\mathbb Z_2$ then the  $\mathbb Z_2$ action 
on $\S^2\times T_0^n$ is given by 
$$t(v,z_1,\ldots,z_n)=(-v,\bar{z}_1,\ldots,\bar{z}_n),\ v\in \S^2, z_j\in \S^1$$ This action is fixed point free and so  there is
 a two-fold covering 
$\S^2\times T^n_0\stackrel{p}{\to}M_n^+-\{x_0\}$ and  a short exact sequence
$$1\to \pi_1(\S^2\times T^n_0)\to\pi_1(M_n^+-\{x_0\})\to\Z_2\to 1$$
Let ${\bf n}$ denote the north pole of $\S^2.$ Then for base points in $\S^2\times T^n_0$ and $M_n^+-\{x_0\}$ we take 
$({\bf n},-1,\ldots,-1)=({\bf n},-\vec{1})$ and $[{\bf n},-1,\ldots,-1]=[{\bf n},-\vec{1}]$ respectively.

This sequence splits. To see this note that the composite 
$\S^2\to \S^2\times T^n_0 \to \S^2$ is  the identity,  where the first map is
$ v\mapsto(v,-\vec{1})$  and the second is just the projection. Both maps are  equivariant with respect to the $\Z_2$-actions, and therefore $ M_n^+-\{x_0\}$ retracts onto $\mathbb R P^2.$

First we consider the case $n=2.$ Choose $-1$ to be the base point in $\S^1.$ The above formula for the action of $\mathbb Z_2$ also defines a $\Z_2$ action on  $\S^2\times(\S^1\vee \S^1).$
This action is fixed point free. The inclusion  $\S^2\times(\S^1\vee \S^1)\subset \S^2\times \S^1\times \S^1$ is equivariant and there exists a $\Z_2$-equivariant strong deformation retract from 
$\S^2\times T^2_0$ onto $\S^2\times(\S^1\vee \S^1).$  Let $a_1,a_2$ be the generators 
$({\bf n},\S^1,-1)$ and $({\bf n},-1,\S^1)$ of  $\pi_1(\S^2\times T^2_0)=\Z*\Z.$ See the Figure below. 

The involution $t:\S^2\times T^2_0\to \S^2\times T^2_0$ induces isomorphisms 
$$\begin{array}{ccc}
\pi_1(\S^2\times(\S^1\vee \S^1),\{{\bf n},-1,-1\})&\stackrel{c}{\to}&
\pi_1(\S^2\times(\S^1\vee \S^1),\{{\bf s},-1,-1\})\\
\pi_1(\S^2\vee(\S^1\vee \S^1),\{{\bf n},-1,-1\})&\stackrel{c}{\to}&
\pi_1(\S^2\vee(\S^1\vee \S^1),\{{\bf n},-1,-1\})
\end{array}$$
where ${\bf s}=-{\bf n}$ is the south pole in $\S^2.$

We have the following commutative diagram
$$\xymatrix{
\S^2\vee_{{\bf n}}(\S^1\vee \S^1)\ar[d]^{i_{{\bf n}}}\ar[rr]^t & & \S^2\vee_{{\bf s}}(\S^1\vee \S^1)\ar[d]^{i_{{\bf s}}} \\
\S^2\times_{{\bf n}}(\S^1\vee \S^1)\ar[rr]^t\ar[rd]_p  & & \S^2\times_{{\bf s}}(\S^1\vee \S^1)\ar[ld]^p  \\
 & M_2^+-\{x_0\} &  }$$
where $i_{{\bf n}}$ and $i_{{\bf s}}$ are inclusions. Here the subscripts ${\bf n}$ and ${\bf s}$ refer to the north and south poles respectively, which we take to be base points of $\S^2$ in the one point unions.  The inclusions $i_{\bf n},i_{\bf s}$  induce isomorphims on $\pi_1$ and therefore 
 $p_*\pi_1(\S^2\vee_{\bf n}(\S^1\vee \S^1))=p_*\pi_1(\S^2\vee_{\bf s}(\S^1\vee \S^1)).$ Thus $t$ sends  $a_1$ to the loop based at $s$ but with
the opposite orientation (similarly for $a_2$). See the Figure below. 

\begin{center}
\includegraphics [scale=0.5] {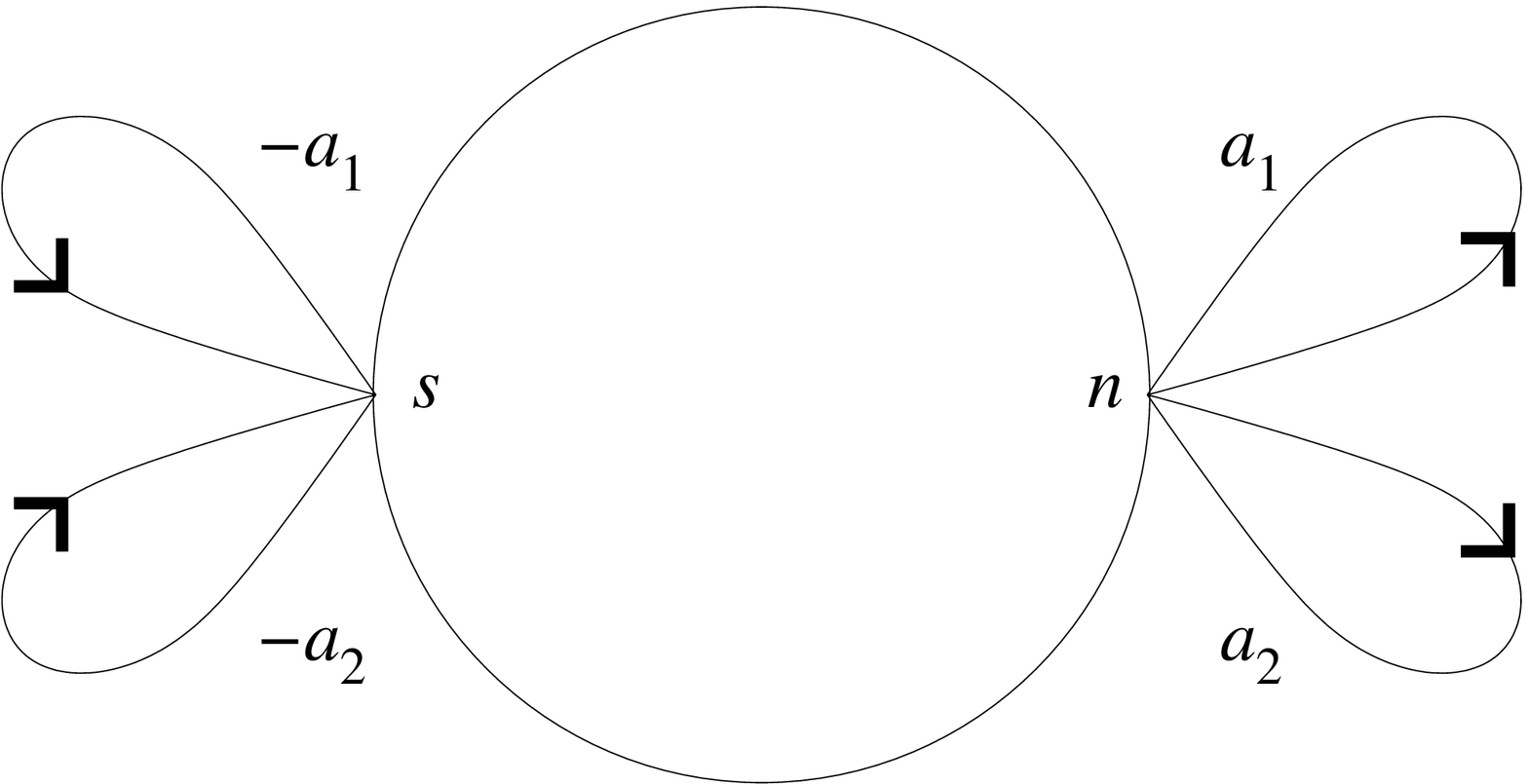}
\end{center}

We now have  $\pi_1(M_2^+-\{x_0\})=
<a_1,a_2,t\ |\ t^2=1, a_1^t=a_1^{-1}, a_2^t=a_2^{-1}>$. 

For the computation of $\pi_1(M_n^+-\{ x_0\}),\ n\ge 3,$ note  that the inclusion $T^n_0 \subset T^n$ induces an isomorphism on $\pi_1.$ Therefore 
$\pi_1(T^n_0)=<a_1,\ldots,a_n\ |\ [a_i,a_j]=1\ \forall\ i,j>.$ 
The various inclusions of $T_0^2$ into $T_0^n$ (corresponding to pairs of generators) show that the action of $t$ on the 
generators is still
given by $a_i^t=a_i^{-1}.$ 
Thus 
$$\pi_1(M_n^+-x_0)=<a_1,\ldots,a_n,t\ |\ t^2=1,[a_i,a_j]=1, a_i^t=a_i^{-1}>,\ 
\mbox{ for $n\geq 3$.}$$
The final step in the calculation of $\pi_1(M_n^+)$ is to use van Kampen's theorem. To do this let 
 $U\subset \S^1$  be a small open connected neighbourhood of $1\in \S^1$ which is invariant under conjugation.  Here small means 
 $-1\not\in U$. Then $N_n=\S^2\times U^n/\sim$ is a contractible neighborhood of $x_0$ in $M_n^+.$ 
We apply van Kampen's theorem to the situation 
$M_n^+=\dis (M_n^+-\{x_0\})\cup N_n.$

The intersection $N_n\cap(M_n^+-\{x_0\})$ is homotopy equivalent to 
$(\S^2\times\S^{n-1})/\Z_2$ where $\Z_2$ acts by multiplication by
$-1$ on both factors. Therefore $\pi_1(N_n\cap(M_n^+-\{x_0\})) \cong\Z$ when $n =2$, and $\Z_2$ when $n\geq 3$.
Thus we need to understand the homomorphism induced by the inclusion 
$N_n\cap(M_n^+-\{ x_0\} )\to M_n^+-\{ x_0\}$.
  
When $n=2$ the inclusion of $N_2\cap(M_2^+-\{x_0\})$ into $M_2^+-\{x_0\}$ induces the following commutative diagram 
$$\xymatrix{
\Z\ar[r]\ar[d]^2 & \Z*\Z\ar[d] \\
\pi_1(N_2\cap(M_2^+-\{x_0\}))\ar[r]\ar[d] & \pi_1(M_2^+-\{x_0\})\ar[d] \\
\Z_2\ar[r]^= & \Z_2 }$$
where the map on top is the commutator map. So if the generator of $\pi_1(N_2\cap(M_2^+-\{x_0\}))$ is sent to 
$w\in\pi_1(M_2^+-\{x_0\})$, then $w^2=[a_1,a_2]$,  and the image of $w$ in $\Z_2$ is $t$. 
Thus we can write $w=a_1^{n_1}a_2^{m_1}\cdots a_1^{n_r}a_2^{m_r}t$ with $n_i,m_i\in\Z$. 
Then 
$$w^2=a_1^{n_1}a_2^{m_1}\cdots a_1^{n_r}a_2^{m_r}a_1^{-n_1}a_2^{-m_1}\cdots a_1^{-n_r}a_2^{-m_r}=
a_1 a_2 a_1^{-1}a_2^{-1}$$
which occurs only if $r=1$ and $n_1=m_1=1$. It follows that $w=a_1a_2 t$.
Thus $$\pi_1(M_2^+)=<a_1,a_2,t\ |\ t^2=1, a_1^t=a_1^{-1}, a_2^t=a_2^{-1}, a_1a_2 t=1>$$ and routine computations 
show that this is the Klein four group.

For $n\geq 3$ the inclusion map $N_n\cap(M_n^+-\{x_0\})\to M_n^+-\{x_0\}$ can be understood by looking at the following diagram
$$\xymatrix{
\S^2\times \S^1\ar[rr]\ar[dd]\ar[rd] & & \S^2\times T^2_0\ar[dd]\ar[rd] & \\
 & \S^2\times S^{n-1}\ar[rr]\ar[dd] & & \S^2\times T^n_0\ar[dd] \\
N_2\cap(M_2^+-\{x_0\})\ar[rr]\ar[rd] & & M_2^+-\{x_0\}\ar[rd] & \\
 & N_n\cap(M_n^+-\{x_0\})\ar[rr] & & M_n^+-\{x_0\}
}$$
Note that the map $N_2\cap(M_2^+-\{x_0\})\to N_n\cap(M_n^+-\{x_0\})$ induces the canonical projection $\Z\to\Z_2$. A chase
argument shows that the inclusion $N_n\cap(M_n^+-\{x_0\})\to M_n^+-\{x_0\}$ imposes the relation $a_1a_2t$ as well, and
therefore 
$$\pi_1(M_n^+)=<a_1,\ldots,a_n,t\ |\ t^2=1,[a_i,a_j]=1, a_i^t=a_i^{-1}, a_1a_2 t=1>.$$
By performing some routine computations we see that this group is isomorphic to $\Z_2^n$. 
This completes the proof of Theorem~\ref{mainth} for $SO(3)$. The cases of $SU(2)$ and $U(2)$ follow from Remark~\ref{cover}. 
\hfill$\Box$

\medskip
Since the map $\pi_1(\vee_n G)\to\pi_1(G^n)$ is an epimorphism, it follows that the inclusion maps 
$$M_n^+(G)\to G^n\hspace{.6cm}if\ \ G=SO(3)$$ 
$$M_n(G)\to G^n\hspace{.5cm}if\ \ G=SU(2),U(2)$$ 
are isomorphisms in $\pi_1$ for all $n\geq 1$. 
Recall that there is a map $Hom(\Gamma,G)\to Map_*(B\Gamma,BG)$, where $Map_*(B\Gamma,BG)$ is the space of pointed 
maps from the classifying space of $\Gamma$ into the classifying space of $G$. Let $Map_*^+(T^n,BG)$ be the component of 
the map induced by the trivial representation. 

\begin{corollary}\label{5.2} The maps
$$M_n^+(G)\to Map_*^+(T^n,BG)\hspace{.5cm}if\ \ G=SO(3)$$ 
$$M_n(G)\to Map_*^+(T^n,BG)\hspace{.5cm}if\ \ G=U(2)$$ 
are injective in $\pi_1$ for all $n\geq 1$.
\end{corollary}
\begin{proof} By induction on $n$, with the case $n=1$ being trivial. Assume $n>1$, and
note that there is a commutative diagram
$$\xymatrix{
M_n^+(SO(3))\ar[r]\ar[d]  & Map_*^+(B\pi_1(T^n),BSO(3))\ar[d] \\
Hom(\pi_1(T^{n-1}\vee \S^1),SO(3))\ar[r]\ar[d]  & Map_*^+(B\pi_1(T^{n-1}\vee \S^1),BSO(3))\ar[d]  \\
Hom(\pi_1(T^{n-1}),SO(3))\times SO(3)\ar[r]  & Map_*^+(B\pi_1(T^{n-1}),BSO(3))\times SO(3) }$$
in which the bottom map is injective in $\pi_1$ by inductive hypothesis, the lower vertical maps are homeomorphisms, 
and the upper left vertical map is injective in $\pi_1$. Thus the map on top is also injective as wanted.
The proof for $U(2)$ is the same. 
\end{proof}

\begin{remark}{\em We have the following observations.
\begin{enumerate}
\item The two-fold cover $\Z_2\to \S^3\times \S^3\to SO(4)$ allows us to study $Hom(\Z^n,SO(4))$. 
Let $M_n^+(SO(4))$ be the component covered by $Hom(\Z^n,\S^3\times \S^3)$. Since  $Hom(\Z^n,\S^3\times \S^3)$ is
homeomorphic to  $Hom(\Z^n,\S^3)\times Hom(\Z^n,\S^3)$, it follows that 
$$\pi_1(M_n^+(SO(4)))=\Z_2^n$$
\item The space $Hom(\Z^2,SO(4))$ has two components. One is $M_2^+(SO(4)),$ which is covered by $\partial^{-1}_{SU(2)^2}(1,1)$,
and the other is covered by $\partial^{-1}_{SU(2)^2}(-1,-1)$, where $\partial_{SU(2)^2}$ is the commutator map of $SU(2)\times SU(2)$.
Recall $\partial^{-1}_{SU(2)}(-1)$ is homeomorphic to $SO(3)$  (see~\cite{AM}), so $\partial^{-1}_{SU(2)^2}(-1,-1)$ is 
homeomorphic to $SO(3)\times SO(3)/\Z_2\times\Z_2$, where the group $\Z_2\times\Z_2$ acts by left diagonal multiplication when it is 
thought of as the subgroup of $SO(3)$ generated by the transformations $(x_1,x_2,x_3)\mapsto(x_1,-x_2,-x_3)$ and
 $(x_1,x_2,x_3)\mapsto(-x_1,x_2,-x_3)$.
\item Corollary~\ref{5.2} holds similarly for $SO(4)$, and trivially for $SU(2)$.

\end{enumerate}
}\end{remark}

\section{\bf Homological Computations}

In this section we compute the $\Z_2$-cohomology of $M_n^+(SO(3))$. The $\Z_2$-cohomology of the other components
of $M_n(SO(3))$ is well-known since these are all homeomorphic to $\S^3/Q_8$.
To perform the computation we will use the description of $M_n^+(SO(3))$ that we
saw in the proof of Theorem~\ref{mainth}. The ingredients we have to consider are the spectral sequence of the fibration 
$\S^2\times T^n_0\to(M_n^+-\{ x_0\})\to \R P^\infty$ whose $E_2$ terms is 
$$\Z_2[u]\otimes E(v)\otimes E(x_1,\ldots,x_n)/(x_1\cdots x_n)$$  
with $deg(u)=(1,0)$, $deg(v)=(0,2)$ and $deg(x_i)=(0,1)$; 
and the spectral sequence of the fibration $\S^2\times\S^{n-1}\to N_n\cap (M_n^+-\{ x_0\})\to\R P^\infty$ whose
$E_2$-term is $$\Z_2[u]\otimes E(v)\otimes E(w)$$
with $deg(u)=(1,0)$, $deg(v)=(0,2)$ and $deg(w)=(0,n-1)$. Note that in both cases $d_2(v)=u^2$, whereas $d_2(x_i)=0$ since 
$H^1(M_n^+-\{ x_0\},\Z_2)=\Z_2^{n+1}$. Therefore the first spectral sequence collapses at the third term.  
As $d_n(w)=u^n$ and $d_j(w)=0$ for $j\neq n$,  the second spectral sequence collapses at the third term when $n=2$
and at the fourth term when $n\geq 3$.  

The last step is to use the Mayer-Vietoris long exact sequence of the pair $(M_n^+-\{ x_0\},N_n)$ which yields the following: for $n=2,3$, 

\bigskip
$$H^q(M_2^+(SO(3)),\Z_2)=\left\{\begin{array}{ccl}
\Z_2 & & q=0\\
\Z_2\oplus\Z_2 & & q=1\\
\Z_2\oplus\Z_2\oplus\Z_2 & & q=2\\
\Z_2\oplus\Z_2\oplus\Z_2 & & q=3\\
\Z_2 & & q=4\\
0 & & q\geq 5 \end{array}\right.$$

\vspace{1cm}
$$H^q(M_3^+(SO(3)),\Z_2)=\left\{\begin{array}{ccl}
\Z_2 & & q=0\\
\Z_2\oplus\Z_2\oplus\Z_2 & & q=1\\
\Z_2\oplus\Z_2\oplus\Z_2\oplus\Z_2\oplus\Z_2\oplus\Z_2 & & q=2\\
\Z_2\oplus\Z_2\oplus\Z_2\oplus\Z_2\oplus\Z_2\oplus\Z_2\oplus\Z_2 & & q=3\\
\Z_2\oplus\Z_2\oplus\Z_2\oplus\Z_2 & & q=4\\
\Z_2 & & q=5\\
0 & & q\geq 6 \end{array}\right.$$

\vspace{1cm}
whereas for $n\geq 4$,

$$H^q(M_n^+(SO(3)),\Z_2)=\left\{\begin{array}{ccl}
\Z_2 & & q=0\\
 & & \\
\Z_2^n & & q=1\\
 & & \\
\Z_2^{{n\choose 1}+{n\choose 2}} & & q=2\\
 & & \\
\dis \Z_2^{{n\choose q-2}+{n\choose q-1}+{n\choose q}} & & 3\leq q\leq n\\
 & & \\
\Z_2^{{n\choose n-1}+1} & & q=n+1\\
 & & \\
\Z_2 & & q=n+2 \\
 & & \\
0 & & q\geq n+3 \end{array}\right.$$

So the Euler characteristic of $M_n^+(SO(3))$ is given by
$$\chi(M_n^+(SO(3)))=\left\{
\begin{array}{ccl}
0 & & n=2\ or\  odd \\
 & & \\
\dis 2+n(n-1)-{n\choose k-1}-{n\choose k}-{n\choose k+1} & & n=2k,\ \ k\geq 2
\end{array}\right.$$

\bigskip
\hfill\
{
\parbox{6cm}{Denis Sjerve\\
{\it Department of Mathematics},\\
University of British Columbia\\
Vancouver, B.C.\\
Canada\\
{\sf sjer@math.ubc.ca}}\
{\hfill}\
\parbox{6cm}{Enrique Torres-Giese\\
{\it Department of Mathematics},\\
University of British Columbia\\
Vancouver, B.C.\\
Canada\\
{\sf enrique@math.ubc.ca}}}\
\hfill

\end{document}